\documentclass[12pt]{article}                                               

\input{amssym.def}                                                           
\input{amssym.tex}                                                           
                                                                             
\begin{document}                                                             
\title{On a  Poisson summation formula for  noncommutative tori}

\author{Igor  ~Nikolaev}


\date{}
 \maketitle


\newtheorem{thm}{Theorem}
\newtheorem{lem}{Lemma}
\newtheorem{dfn}{Definition}
\newtheorem{rmk}{Remark}
\newtheorem{cor}{Corollary}
\newtheorem{prp}{Proposition}
\newtheorem{exm}{Example}
\newtheorem{cnj}{Conjecture}

\newcommand{\ch}{\hbox{\bf Char}}
\newcommand{\sa}{\hbox{\bf sa}}
\newcommand{\ka}{\hbox{\bf k}}
\newcommand{\n}{\hbox{\bf n}}
\newcommand{\mod}{\hbox{\bf mod}}
\newcommand{\Z}{\hbox{\bf Z}}
\newcommand{\Q}{\hbox{\bf Q}}
\newcommand{\tr}{\hbox{\bf tr}}

\begin{abstract}
It is proved that a maximal abelian subalgebra of the noncommutative torus
commutes with the Laplace operator on a complex torus.  As a corollary,
one gets an analog of the Poisson summation formula for noncommutative tori.

\vspace{7mm}

{\it Key words and phrases: Selberg trace formula,  noncommutative  torus}

\vspace{5mm}
{\it MSC: 11F72 (Selberg trace formula);  46L85 (noncommutative topology)}
\end{abstract}

\section{Introduction}
The Poisson summation formula is an elementary and fundamental fact of harmonic analysis  
and representation theory. The simplest case of such a formula says that for every
function $f\in C_0^{\infty}({\Bbb R})$  it holds
\begin{equation}\label{eq1}
\sum_{n\in {\Bbb Z}}f(n)=\sum_{n\in {\Bbb Z}}\hat f(n),
\end{equation}
where $C_0^{\infty}({\Bbb R})$ is the set of $C^{\infty}$-smooth  complex-valued functions 
with a compact support  on the real line ${\Bbb R}$ and
$\hat f(\nu):=\int_{-\infty}^{\infty}f(x)e^{-2\pi i\nu x}dx$ is  the Fourier transform of  $f$. 
Let ${\Bbb H}=\{x+iy\in {\Bbb C}~|~y>0\}$ be the Lobachevsky half-plane and $\tau\in {\Bbb H}$;
an analog of formula (\ref{eq1}) for the two-dimensional lattice  $L_{\tau}={\Bbb Z}+\tau {\Bbb Z}$
can be written as
\begin{equation}\label{eq3}
\sum_{m\in {\Bbb Z}}\sum_{n\in {\Bbb Z}}f(m^2+2mn~\Re (\tau)+n^2|\tau|^2)
=\sum_{r=0}^{\infty}\lambda(r)\hat f(r),
\end{equation}
where function $f\in C_0^{\infty}({\Bbb R}^2)$ is radially symmetric
(i.e. $f(u,v)=Const$  on  $u^2+v^2=r\ge 0$)  while   $\lambda(r):=|\{(m,n)\in {\Bbb Z}^2 : m^2+n^2=r\}|$
the  multiplicity function and
$\hat f$ the Fourier transform of $f$ given by the formula
\begin{equation}\label{eq4}
\hat f(r)={\pi\over\Im (\tau)}\int_0^{\infty}f(s)
J_0\left[2\pi\sqrt{s}\sqrt{{m^2|\tau|^2-2mn~\Re (\tau)+n^2\over \Im^2 (\tau)}}\right]ds, 
\end{equation}
with  $J_0(z)={1\over\pi}\int_0^{\pi}\cos(z\cos\alpha)d\alpha$ being the Bessel function.

Recall that each radially symmetric function  $f\in C_0^{\infty}({\Bbb R}^2)$ gives rise
to a symmetric Hilbert-Schmidt integral operator on the Hilbert space $L^2({\Bbb C}/L_{\tau})$
acting by the formula: 
\begin{equation}\label{eq4bis}
(T_f\varphi)(z)=\int_{{\Bbb C}/L_{\tau}} f(z,w)\varphi(w)dw,
\end{equation}
 where $f(z,w)=\sum_{z_0\in L_{\tau}}f(z+z_0, w)$,  see e.g. [Iwaniec 1995]  \cite{I}, p. 5.
 Denote by $\sum f(L_{\tau})$ the LHS of (\ref{eq3});   a link between $T_f$ and the Poisson summation
 is given by the formula
\begin{equation}\label{eq4bisbis}
\tr ~(T_f)=\sum f(L_{\tau}),
\end{equation}
 where $\tr$ is the trace of $T_f$.   
 The $T_f$  commutes with the Laplace operator 
 $\Delta={\partial^2\over\partial x^2}+{\partial^2\over\partial y^2}$
 on the  complex torus ${\Bbb C}/L_{\tau}$ and operators $T_f$ commute with each other 
 for all  $f\in C_0^{\infty}({\Bbb R})$, {\it ibid.}    It is easy to see,  that complex conjugation 
 defines an adjoint operator $T^*_f =T_{\bar f}$.   The norm closure of the $\ast$-algebra
 generated by all $T_f$ is a commutative $C^*$-algebra,   see e.g.  [Murphy 1990] \cite{M}
 for an introduction;   such a $C^*$-algebra we denote by 
 \begin{equation}
 {\cal R}({\Bbb C}/L_{\tau}):=\overline{\{T_f : f\in C_0^{\infty}({\Bbb R}^2)\}}.
 \end{equation}

In this note we construct an inclusion  of the algebra  ${\cal R}({\Bbb C}/L_{\tau})$
into  a  noncommutative torus ${\cal A}_{\theta}$,
i.e. the  $C^*$-algebra  generated by  unitary operators  $u$ and $v$ satisfying the commutation 
relation $vu=e^{2\pi i\theta}uv$ for a constant $\theta\in {\Bbb R}$  [Rieffel 1990]  \cite{Rie2}. 
 Namely, let $[a_0, a_1,  \dots]$ be the regular 
 continued fraction of $\theta$ and consider the Bratteli diagram in Figure 1, 
where $a_i$ is the multiplicity of edges of the diagram  [Bratteli  1972]  \cite{Bra1}.  
(Note that the diagram is infinite, unless $\theta$ is a  rational number;  in this case 
the canonical trace on ${\cal A}_{\theta}$ is defined by continuity from the irrational values 
of $\theta$.)   Let $X_{\theta}$ be the Bratteli compactum,  i.e. a Cantor set obtained from 
 the infinite paths of the Bratteli diagram,   see [Herman, Putnam \& Skau  1992]  \cite{HePuSk1}, pp. 837-838. 
We shall denote by ${1\over\mu} X_{\theta}$ the Cantor set $X_{\theta}$ endowed with
the measure $\mu$ and by $C({1\over\mu} X_{\theta})$  the commutative $C^*$-algebra of complex-valued 
functions on ${1\over\mu} X_{\theta}$.  We shall write  $({\cal A}_{\theta}, {1\over\mu} e)=F({\Bbb C}/L_{\tau})$
to denote the noncommutative torus ${\cal A}_{\theta}$ with a scaled unit  ${1\over\mu} e$ corresponding to the 
complex torus ${\Bbb C}/L_{\tau}$ under a functor $F$, see Section 2.1. 
 Our main results  can be stated as follows.

\begin{figure}
\begin{picture}(300,60)(-40,0)
\put(110,30){\circle{3}}
\put(120,20){\circle{3}}
\put(140,20){\circle{3}}
\put(160,20){\circle{3}}
\put(120,40){\circle{3}}
\put(140,40){\circle{3}}
\put(160,40){\circle{3}}

\put(110,30){\line(1,1){10}}
\put(110,30){\line(1,-1){10}}
\put(120,42){\line(1,0){20}}
\put(120,40){\line(1,0){20}}
\put(120,38){\line(1,0){20}}
\put(120,40){\line(1,-1){20}}
\put(120,20){\line(1,1){20}}
\put(140,41){\line(1,0){20}}
\put(140,39){\line(1,0){20}}
\put(140,40){\line(1,-1){20}}
\put(140,20){\line(1,1){20}}

\put(180,20){$\dots$}
\put(180,40){$\dots$}

\put(125,52){$a_0$}
\put(145,52){$a_1$}

\end{picture}

\caption{Bratteli diagram of ${\cal A}_{\theta}$.}
\end{figure}
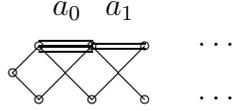

\noindent
 \begin{thm}\label{thm1}
 There exists a trace-preserving isomorphism: 
\begin{equation}
{\cal R}({\Bbb C}/L_{\tau})\cong C\left({1\over\mu} X_{\theta}\right)
\subset  \left({\cal A}_{\theta}, {1\over\mu} e\right),
\end{equation}
where
 $C({1\over\mu} X_{\theta})$ is the maximal abelian  subalgebra
 of the $C^*$-algebra $({\cal A}_{\theta}, {1\over\mu} e)$.   
  \end{thm}
\begin{cor}\label{cor1}
{\bf (Poisson summation formula)}
For each radially symmetric function $f\in C_0^{\infty}({\Bbb R}^2)$  
 there exists a Hilbert-Schmidt  operator $T_ f\in ({\cal A}_{\theta}, {1\over\mu} e)$, 
 such that:  
\begin{equation}
\sum f(L_{\tau})=\tr~(T_f),
\end{equation}
 where $\tr$  is the  canonical trace on  $({\cal A}_{\theta}, {1\over\mu} e)$ 
 [Rieffel  1981]  \cite{Rie1}.   The operator $T_f$ is self-adjoint if and
 only if $f$ is a real-valued function. 
 \end{cor}
\begin{rmk}\label{rmk1}
{\normalfont
Theorem \ref{thm1} says that  the $C^*$-algebra
${\cal R}({\Bbb C}/L_{\tau})$  is  a  maximal abelian subalgebra (masa) of the 
noncommutative torus  $({\cal A}_{\theta}, {1\over\mu} e)=F({\Bbb C}/L_{\tau})$;
 the masa determines the $C^*$-algebra $({\cal A}_{\theta}, {1\over\mu} e)$
itself  by taking  the crossed product  of $C({1\over\mu} X_{\theta})$
with the Vershik homeomorphism of the Cantor set 
${1\over\mu} X_{\theta}$, see Section 2.2.
}
\end{rmk}
 The structure of the article is as follows.  We recall some useful facts in
 Section 2.  Theorem \ref{thm1} and corollary \ref{cor1} are proved in Section 3.

\section{Preliminaries}
We  briefly review a relation between complex and noncommutative tori, the $C^*$-dynamical systems
on the Cantor set  and the Selberg trace formula.  For an extended account and details we refer the reader 
to \cite{Nik1},  [Herman, Putnam \& Skau  1992]  \cite{HePuSk1} and    [Iwaniec 1995]  \cite{I}, respectively.

\subsection{Complex and noncommutative tori}
Recall that   ${\Bbb C}/L_{\tau}$ is isomorphic  to the intersection two quadric 
surfaces in the complex projective space of the form
 $\{(u,v,w,z)\in {\Bbb C}P^3 ~|~u^2+v^2+w^2+z^2 ={1-\alpha\over 1+\beta}v^2+
 {1+\alpha\over 1-\gamma}w^2+z^2=0\}$,  where $\alpha,\beta,\gamma$ are some complex constants 
 such that $\alpha+\beta+\gamma+\alpha\beta\gamma=0$. 
It was proved by Sklyanin that a free ${\Bbb C}$-algebra $S_{\alpha,\beta,\gamma}$  on four generators $x_i$ and   
six quadratic relations
\begin{equation}
\left\{
\begin{array}{ccc}
x_1x_2-x_2x_1 &=& \alpha(x_3x_4+x_4x_3),\\
x_1x_2+x_2x_1 &=& x_3x_4-x_4x_3,\\
x_1x_3-x_3x_1 &=& \beta(x_4x_2+x_2x_4),\\
x_1x_3+x_3x_1 &=& x_4x_2-x_2x_4,\\
x_1x_4-x_4x_1 &=& \gamma(x_2x_3+x_3x_2),\\ 
x_1x_4+x_4x_1 &=& x_2x_3-x_3x_2,
\end{array}
\right.
\end{equation}
is a {\it coordinate ring}   of complex torus ${\Bbb C}/L_{\tau}$,  i.e.  
$\hbox{{\bf Mod}}~(S_{\alpha,\beta,\gamma})/
\hbox{{\bf Tors}}\cong \hbox{{\bf Coh}}~({\Bbb C}/L_{\tau})$,
 where {\bf Coh} is  the category of quasi-coherent sheaves on ${\Bbb C}/L_{\tau}$, 
  {\bf Mod}  the category of graded left modules over the graded ring $S_{\alpha,\beta,\gamma}$
 and  {\bf Tors}  the full sub-category of {\bf Mod} consisting of the
torsion modules,  see e.g.  [Smith \& Stafford 1993]  \cite{SmiSta1}, p.267.    
The closure of a self-adjoint representation of the Sklyanin algebra $S_{\alpha,\beta,\gamma}$
by linear operators on a Hilbert space ${\cal H}$ is  isomorphic to the algebra  ${\cal A}_{\theta}$ with a scaled unit ${1\over\mu}e$ for a constant 
$\mu>0$.  The bijection  $F: {\Bbb C}/L_{\tau}\mapsto ({\cal A}_{\theta}, {1\over\mu}e)$ is a functor 
from  isomorphic complex  tori to the  stably isomorphic noncommutative tori 
\cite[Section 1.3]{Nik1}.

\subsection{$C^*$-dynamical systems on the Cantor set}
An {\it $AF$-algebra}  (approximately finite-dimensional  $C^*$-algebra) is defined to
be the  norm closure of an ascending sequence of the finite-dimensional
$C^*$-algebras $M_n$'s, where  $M_n$ is the $C^*$-algebra of the $n\times n$ matrices
with the entries in ${\Bbb C}$.
To describe the ascending sequence, we use 
an infinite graph called a {\it Bratteli diagram} of the
$AF$-algebra   [Bratteli  1972]    \cite{Bra1}.    The Bratteli diagram defines a unique  $AF$-algebra.
By an ordered Bratteli diagram  one understands a natural (partial) order between the 
edges of the diagram defined as follows.  Two edges $e$ and $e'$ are comparable,
if and only if, there exist two paths through $e$ and $e'$ respectively, which have an edge
in common;  in this case one writes $e>e'$ whenever $e$ lies above $e'$ on the diagram;
we refer the reader to   [Herman, Putnam \& Skau  1992]  \cite{HePuSk1}, pp. 835-836 for
the details.

To each Bratteli diagram one assigns a Cantor set, i.e. the compact totally disconnected
metric space $X$ consisting of the infinite paths of the diagram.   Let $P_{0,k}$ denote
a path $(e_0,\dots, e_k)$ consisting of $k+1$ edges $e_i$ of the Bratteli diagram.  
The inverse limit $X=\lim P_{0,k}$ is a Cantor set defined by the discrete topology
on $P_{0,k}$,  see [Herman, Putnam \& Skau  1992]  \cite{HePuSk1}, pp. 837-838.  
The Cantor set $X$ is called a {\it Bratteli compactum}.

The ordered Bratteli diagram gives rise to a homeomorphism $\varphi: X\to X$
of the Bratteli compactum $X$ associated to the diagram.  Roughly speaking,
$\varphi$ sends an infinite path of the Bratteli diagram (i.e. a point of $X$)
to the same infinite path but with a finite number of edges replaced by the 
successor edges defined by the ordering,  see [Herman, Putnam \& Skau  1992]  \cite{HePuSk1},
p. 838.  The map $\varphi$ is called a {\it Vershik homomorphism} of the Bratteli
compactum $X$. The iterations of $\varphi$ define a minimal dynamical system
$(X,\varphi)$ on the Cantor set $X$.  By a  $C^*$-dynamical system on the 
Cantor set $X$ one understands the crossed product $C^*$-algebra 
$C(X)\rtimes_{\varphi} {\Bbb Z}$,  where $C(X)$ is the $C^*$-algebra of 
complex-valued functions on $X$.    (This is true also for the rational values of $\theta$.)  
\begin{lem}\label{lm2}
$C(X)\rtimes_{\varphi} {\Bbb Z}\subset {\Bbb A}$, where ${\Bbb A}$ 
is the  $AF$-algebra defined by  the Bratteli diagram upon which
the Cantor set $X$ was constructed.  
\end{lem}
{\it Proof.}  See [Putnam 1989]  \cite{Put1}, pp. 346-350. 
$\square$

\subsection{Selberg trace formula}
To put the Poisson summation formula in a general context,  consider the
homogeneous space $G/K$, where $G\cong SL(2, {\Bbb R})$ ($G\cong {\Bbb R}^2$,  resp.)
and $K\cong SO(2, {\Bbb R})$ ($K$ is trivial, resp.) are the Lie groups; 
thus $G/K\cong {\Bbb H}$  is the Lobachevsky half-plane
($G/K\cong {\Bbb R}^2$ is the Euclidean plane, resp.) 
Recall that the left regular representation of $G$ gives rise to a linear operator
$(T_{\gamma} f)(z)=f(\gamma z)$ on the space of complex-valued functions
$f: G\to {\Bbb C}$.  A linear operator $T$ is said to be {\it invariant} if it commutes
with $T_{\gamma}$, i.e. $T(f(\gamma z))=(Tf)(\gamma z), \quad \forall\gamma\in G$.
\begin{exm}
\textnormal{
The Laplace-Beltrami operator $\Delta=y^2
\left({\partial^2\over\partial x^2}+{\partial^2\over\partial y^2}\right)$
(the Laplace operator $\Delta={\partial^2\over\partial x^2}+{\partial^2\over\partial y^2}$,
resp.)  is an invariant differential operator on $C^{\infty}(G/K)$. 
}
\end{exm}
\begin{exm}
\textnormal{
Consider the integral operator  $(T_f\varphi)(z)=\int_{G/K} f(z,w)\varphi(w)dw$,
where the kernel $f: G/K\times G/K\to {\Bbb C}$ satisfies condition 
$f(\gamma z, \gamma w)=f(z,w)$ for all $\gamma\in G$;    the
$T_f$ is an invariant integral operator and such operators  commute with each other.  
}
\end{exm}
\begin{lem}\label{lm1}
The invariant integral operators $T_f$ commute with the Laplace-Beltrami
operator $\Delta$ (the Laplace operator $\Delta$, resp.);  moreover,
any eigenfunction of $\Delta$ coincides with an eigenfunction of all
invariant integral operators.
\end{lem}
{\it Proof.}  See e.g.  [Iwaniec 1995]  \cite{I}, pp. 29-33.
$\square$

\bigskip
Let $\Gamma$ be a discrete subgroup of $SL(2, {\Bbb R})$
(a lattice $L_{\tau}={\Bbb Z}+{\Bbb Z}\tau$ of ${\Bbb R}^2\cong {\Bbb C}$, resp.);
then $\Gamma\backslash {\Bbb H}$ is a Riemann surface
(${\Bbb C}/L_{\tau}$ is a complex torus, resp.)
A function $f: G/K\to {\Bbb C}$ is said to be {\it automorphic}
with respect to $\Gamma$ if it satisfies the periodicity condition
$f(\gamma z)=f(z)$ for all $\gamma\in\Gamma$;  we shall denote
the space of such functions by $A(G/K)$.  Each function 
$f: G/K\to {\Bbb C}$ with a rapid decay on $G/K$ (e.g. a radially symmetric function $f\in C^{\infty}_0(G/K)$)
gives rise to an automorphic function with respect to $\Gamma$ given
by the formula $f(z):=\sum_{\gamma\in\Gamma}f(\gamma z)$.

Further we  restrict  the Laplace operator $\Delta$ to  the space $A(G/K)$
and  consider the invariant integral operators $T_f$ on $A(G/K)$  given by the formula
\begin{equation}
(T_f\varphi)(z)=\int_{\Gamma\backslash G/K}  f(z,w) \varphi(w)dw,
\end{equation}
where $f(z,w)=\sum_{\gamma\in\Gamma} f(z^{-1}\gamma w)$;
it is known that $T_f$ are compact trace class operators and lemma \ref{lm1}
is valid for them.  Moreover,  the operators $T_f$ satisfy the 
following remarkable formula.
\begin{lem}\label{sel}
{\bf (Selberg trace formula)}
If $T_f$ is a  self-adjoint operator,  then
\begin{equation}\label{selberg}
\tr~(T_f)=\sum_{[\gamma]\in\Gamma} \int_{F_{\gamma}} f([\gamma])dz 
\end{equation}
where $[\gamma]:=\{\tau^{-1}\gamma\tau ~|~ \tau\in\Gamma\}$  is the conjugacy class of $\gamma$
and $F_{\gamma}$ is the fundamental domain of the centralizer of $\gamma$.
\end{lem}
{\it Proof.}  See e.g.  [Iwaniec 1995]  \cite{I}.
$\square$

\begin{rmk}\label{rmk2}
{\normalfont
The Selberg trace formula for $G\cong {\Bbb C}$ and $\Gamma=L_{\tau}$
coincides with the Poisson summation formula.  Indeed,    the LHS of (\ref{selberg}) 
is equal to the $\sum \hat f(L_{\tau})$,  where $\hat f$ is the Fourier 
transform of $f$.
The RHS of (\ref{selberg}) corresponds to  the sum $\sum f(L_{\tau})$. 
}
\end{rmk}

\section{Proofs}
\subsection{Proof of theorem \ref{thm1}}
We shall split the proof in a series of lemmas.
\begin{lem}\label{lm3}
${\cal R}({\Bbb C}/L_{\tau})$ is commutative $C^*$-subalgebra of an $AF$-algebra ${\Bbb A}$.
\end{lem}
{\it Proof.} Let ${\cal R}_0({\Bbb C}/L_{\tau})$ be a dense $\ast$-subalgebra 
of the $C^*$-algebra ${\cal R}({\Bbb C}/L_{\tau})$ consisting of the Hilbert-Schmidt
integral operators $T_f$,  where  $f(r)\in C_0^{\infty}({\Bbb R})$ is a radially symmetric function.   

Denote by $V_i\subset L^2({\Bbb C}/L_{\tau})$ the eigenspace of a Hilbert-Schmidt 
operator $T_f\in {\cal R}_0({\Bbb C}/L_{\tau})$ corresponding to an eigenvalue $\lambda_i$.  
It is well known that $dim~V_i:=n_i<\infty$ and $V_i$ is a common eigenspace 
for all operators $T_f\in {\cal R}_0({\Bbb C}/L_{\tau})$,  see e.g.  [Iwaniec 1995] \cite{I},  p.210.
(The $V_i$ is also an eigenspace for the Laplace operator $\Delta$.)

A restriction of the algebra   ${\cal R}_0({\Bbb C}/L_{\tau})$ to the subspace $V_i$
is a commutative subalgebra of the matrix algebra $M_{n_i}({\Bbb C})$.  

Denote by $\lambda_1>\lambda_2>\dots$ the sequence of eigenvalues 
of the Hilbert-Schmidt operator $T_f$, see  [Iwaniec 1995] \cite{I},  p.210;  
let $\{V_1, V_2,\dots\}$ be the corresponding  eigenspaces of dimensions 
$\{n_1, n_2,\dots\}$,  respectively.  Consider the finite-dimensional $C^*$-algebra
\begin{equation}
{\Bbb A}_i=M_{n_1}({\Bbb C})\oplus\dots\oplus M_{n_i}({\Bbb C}).
\end{equation}
Similarly,  the restriction of  ${\cal R}_0({\Bbb C}/L_{\tau})$ to the sum
of vector spaces $V_1\oplus\dots\oplus V_i\subset L^2({\Bbb C}/L_{\tau})$ 
is a commutative subalgebra of the $C^*$-algebra ${\Bbb A}_i$.

One gets an infinite ascending sequence of the finite-dimensional $C^*$-algebras
${\Bbb A}_1\subset {\Bbb A}_2\subset\dots$;   the corresponding $AF$-algebra we 
shall denote by ${\Bbb A}$.

Because $\cup {\Bbb A}_i$ is a dense subalgebra of ${\Bbb A}$,   the closure of 
the $\ast$-algebra  ${\cal R}_0({\Bbb C}/L_{\tau})$ defines a commutative $C^*$-subalgebra
of the $AF$-algebra ${\Bbb A}$;   in other words,    the commutative 
$C^*$-algebra ${\cal R}({\Bbb C}/L_{\tau})$ is a subalgebra  of the $AF$-algebra ${\Bbb A}$.  
Lemma \ref{lm3} follows.
 $\square$

\begin{lem}\label{lm3.5}
The algebra ${\cal R}({\Bbb C}/L_{\tau})\cong C(X)$,  where $X$ is the Bratteli compactum
of algebra ${\Bbb A}$. 
\end{lem}
{\it Proof.} 
Because ${\cal R}({\Bbb C}/L_{\tau})$ is a commutative  $C^*$-algebra,
 one gets by the Gelfand theorem that ${\cal R}({\Bbb C}/L_{\tau})\cong C(X_0)$,
 where $X_0$ is a compact Hausdorff topological space.  
 
 On the other  hand,  we showed earlier that the algebra  ${\cal R}({\Bbb C}/L_{\tau})$
 defines the $AF$-algebra ${\Bbb A}$, see lemma \ref{lm3}.
 
 The only commutative $C^*$-subalgebra of ${\Bbb A}$ with such a property 
 is the $C^*$-algebra $C(X)$,  where $X$ is the Bratteli compactum of ${\Bbb A}$,
 see [Herman, Putnam \& Skau  1992]  \cite{HePuSk1},  Theorem 8.8.  

Therefore $X=X_0$ and  lemma \ref{lm3.5} follows. 
$\square$

\begin{rmk}\label{rmk2.5}
{\normalfont
${\cal R}({\Bbb C}/L_{\tau})$ is an $AF$-algebra. 
}
\end{rmk}
{\it Proof.}  In view of lemma \ref{lm3.5},  we have  ${\cal R}({\Bbb C}/L_{\tau})\cong C(X)$,  
where $X$ is a Cantor set.   It is known, that in this case ${\cal R}({\Bbb C}/L_{\tau})$
is an $AF$-algebra,  see e.g. [Herman, Putnam \& Skau  1992]  \cite{HePuSk1}, 
p. 831.  Remark \ref{rmk2.5} follows. 
$\square$

\begin{lem}\label{lm4}
${\cal R}({\Bbb C}/L_{\tau})\cong C({1\over\mu} X_{\theta})$. 
\end{lem}
{\it Proof.} 
Consider a crossed product 
\begin{equation}\label{123}
{\cal R}({\Bbb C}/L_{\tau})\rtimes {\Bbb Z}\subset {\Bbb A}
\end{equation}
by the Veshik homeomorphism of the Bratteli compactum $X$ of the
$AF$-algebra ${\Bbb A}$;  the crossed product is not an $AF$-algebra 
but its ordered $K_0$-group is isomorphic to such of ${\Bbb A}$,  see 
[Herman, Putnam \& Skau  1992]  \cite{HePuSk1}.

On the other hand,  it is known that the $AF$-algebra is defined up to an isomorphism 
by its ordered $K_0$-group,  see [Elliott 1976]  \cite{Ell1}.  The crossed product (\ref{123})
has  a unique canonical trace;  therefore  one gets a dense abelian subgroup of the real line 
of the form: 
\begin{equation}\label{1234}
K_0({\cal R}({\Bbb C}/L_{\tau})\rtimes {\Bbb Z})\cong {\Bbb Z}\lambda_1+\dots+{\Bbb Z}\lambda_k,
\end{equation}
where $\lambda_i\in {\Bbb R}$ are some constants and the positive cone  
is defined by the inequality  ${\Bbb Z}\lambda_1+\dots+{\Bbb Z}\lambda_k >0$.  
But the algebra ${\cal R}({\Bbb C}/L_{\tau})$ depends on the complex modulus $\tau=x+iy\in {\Bbb H}$ whose real 
dimension is $2$;  therefore in formula (\ref{1234}) we have $k=2$,  i.e. only $\lambda_1$ and $\lambda_2$
are independent parameters. 

Put $\mu=\lambda_1>0$ and $\theta={\lambda_2\over\lambda_1}$.  Then 
\begin{equation}\label{12345}
K_0({\cal R}({\Bbb C}/L_{\tau})\rtimes {\Bbb Z})\cong \mu({\Bbb Z}+{\Bbb Z}\theta)\cong K_0(({\Bbb A}_{\theta}, {1\over\mu}e)),
\end{equation}
where ${\Bbb A}_{\theta}$ is the $AF$-algebra defined by the Bratteli diagram in Figure  1.  
But  $K_0(({\Bbb A}_{\theta}, {1\over\mu}e))\cong K_0(C({1\over\mu} X_{\theta}) \rtimes {\Bbb Z})$
and comparing with (\ref{12345}) one gets a trace-preserving  isomorphism ${\cal R}({\Bbb C}/L_{\tau})\cong C({1\over\mu} X_{\theta})$
of the corresponding masas;  this   fact  follows from 
the strong orbit equivalence of Cantor minimal systems established in [Herman, Putnam \& Skau 1992]  \cite{HePuSk1} 
and  Elliott's description of the Choquet simplex of the 
tracial states on the $AF$-algebras [Elliott 1976] \cite{Ell1}. 
Lemma \ref{lm4} follows. 
$\square$

\begin{lem}\label{lm5}
The $C({1\over\mu} X_{\theta})$ is a maximal abelian subalgebra of the $C^*$-algebra  
$({\cal A}_{\theta}, {1\over\mu}e)$.
\end{lem}
{\it Proof.}  
  Recall that 
\begin{equation}
C\left({1\over\mu} X_{\theta}\right) \rtimes_{\varphi} {\Bbb Z}
\subset \left({\Bbb A}_{\theta},  {1\over\mu}  e\right)
\subset \left({\cal A}_{\theta}, {1\over\mu} e\right),
\end{equation}
where $\varphi$ is the Vershik homeomorphism, 
see [Putnam 1989]  \cite{Put1},  Theorem 6.7 
for the first inclusion and [Pimsner \& Voiculescu 1980]  \cite{PiVo1} for the second inclusion. 
On the other hand,   the maximal abelian subalgebra of 
$C({1\over\mu } X_{\theta})\rtimes_{\varphi} {\Bbb Z}$
is isomorphic to the $C^*$-algebra  $C({1\over\mu} X_{\theta})$. 
Thus one  gets the following inclusions  
 \begin{equation}\label{eq17}
 C\left({1\over\mu} X_{\theta}\right)\subset
 C\left({1\over\mu } X_{\theta}\right)\rtimes_{\varphi} {\Bbb Z}\subset
\left({\Bbb A}_{\theta},  {1\over\mu} e\right)\subset 
\left({\cal A}_{\theta}, {1\over\mu}e\right).
\end{equation}
On the other hand,  it is known that $K_0(X_{\theta}\rtimes_{\varphi} {\Bbb Z})\cong
K_0({\Bbb A}_{\theta})\cong  K_0({\cal A}_{\theta})$. 
Therefore it follows from (\ref{eq17}) that $C({1\over\mu} X_{\theta})$ is a maximal abelian subalgebra 
of the larger 
algebra  $({\cal A}_{\theta}, {1\over\mu}e)$.  
 Lemma \ref{lm5} follows.
 $\square$

 \bigskip
 Theorem \ref{thm1} follows from lemmas \ref{lm4} and \ref{lm5}.
 $\square$

\subsection{Proof of corollary \ref{cor1}}
Let $f\in C^{\infty}_0({\Bbb R})$ be a radially symmetric function 
on the lattice $L_{\tau}$.  Since the function is compactly supported, 
it defines a Hilbert-Schmidt integral  operator $T_f$ on the Hilbert 
space $L^2({\Bbb C}/L_{\tau})$.    In view of (\ref{eq4bisbis}),  one gets
\begin{equation}
\sum f(L_{\tau})=\tr~(T_f). 
\end{equation}
On the other hand,  theorem \ref{thm1} says that  $T_f \in  ({\cal A}_{\theta}, {1\over\mu}e)$.
It remains to show that $\tr~(T_f)$ coincides with the canonical trace on the noncommutative 
torus. 

Indeed, take a representation of $T_f$ on the Hilbert space $L^2({\Bbb C}/L_{\tau})$
such that the algebra ${\cal R}({\Bbb C}/L_{\tau})$ has a unique trace. 
Then the noncommutative torus $({\cal A}_{\theta}, {1\over\mu}e)$ admits a representation
on the same space. These two representation are different in general, but  
 always unitarily equivalent. In particular, the value of $\tr~(T_f)$ is preserved for any unitary transformation
of the space $L^2({\Bbb C}/L_{\tau})$. 
Corollary \ref{cor1} follows.  
$\square$



\vskip1cm
\textsc{Department of Mathematics and Computer Science, St.~John's University, 8000 Utopia Parkway,  
New York,  NY 11439, United States;} ~\textsc{E-mail:} {\sf igor.v.nikolaev@gmail.com}

\end{document}